\documentclass[12pt]{article}
\usepackage{graphicx}
\usepackage{amsthm,amssymb}

{\newtheorem{definition}{Definition}
\newtheorem{lemma}{Lemma}
\newtheorem{theorem}{Theorem}

}

{\theoremstyle{definition}
\newtheorem{example}{Example}
\newtheorem*{remark}{Remark}
\newtheorem*{conventions}{Conventions}}

\def\BeginExample{\begin{example}}
\def\EndExample{\end{example}}
\def\BeginRemark{\begin{remark}}
\def\EndRemark{\end{remark}}

\def\BeginProofOf#1{\noindent{\bf Proof of #1}: }
\def\EndProof{\QED\betweenskip}
\def\betweenskip{\vskip10pt}

\def\boxitat#1#2{\vbox             
     {\hrule\hbox{\vrule\kern#1%
     \vbox{\kern #1\hbox{#2}\kern#1}\kern#1\vrule}\hrule}}
\def\enclose#1{\boxitat{0pt}{#1}}

\def\qedspace{\null}
\def\qedblack{\qedspace                    
  \lower.6pt\hbox{\vrule height7pt width 5pt}}
\def\qedwhite{\qedspace                    
  \lower.6pt\enclose{%
        \hbox{\vrule height7pt width 0pt\hskip6pt}}}
\def\QED{\qedwhite}

\def\beq{\begin{equation}}
\def\eeq{\end{equation}}

\def\naturals{{\mathbb N}}
\def\integers{{\mathbb Z}}
\def\partsum{\Sigma}
\def\len{{\rm len}}

\def\lrf{\phi}  
\def\lrfs{\Phi} 

\def\set#1{{\cal #1}}
\def\mod{{\rm \,mod\,}}
\def\p{u}

\def\cvec#1{{\vec{\mathbf{#1}}}}
\def\Tr{{\rm t}}
\def\free{\varphi}

\def\sC{{\cal C}}  

\def\NC{C} 

\def\L{{\mathbf L}}  
\def\bfR{{\mathbf R}} 
\def\calM{{\cal M}}  

\def\|{\parallel}

\def\Ex{{\mathsf E}}

\def\km{{\lfloor \log m-\omega(m) \rfloor}}
\def\kp{{\lfloor \log m+\omega(m) \rfloor}}
\def\kmn{{\lfloor \log n-\omega(n) \rfloor}}
\def\kpn{{\lfloor \log n+\omega(n) \rfloor}}

\title {Locally Restricted Compositions IV.\\
Nearly Free Large Parts and Gap-Freeness}

\author{Edward A. Bender\\
\small Department of Mathematics\\[-0.8ex]
\small University of California, San Diego\\[-0.8ex]
\small La Jolla, CA 92093-0112\\
\small\tt ebender@ucsd.edu\\
\and
E. Rodney Canfield\thanks
{Research supported by NSA Mathematical Sciences Program.}\\
\small Department of Computer Science\\[-0.8ex]
\small University of Georgia\\[-0.8ex]
\small Athens, GA 30602\\
\small\tt erc@cs.uga.edu\\
\and
Zhicheng Gao\thanks
{Research supported by NSERC.}\\
\small School of Mathematics and Statistics\\[-0.8ex]
\small Carleton University\\[-0.8ex]
\small Ottawa, Ontario K1S5B6\\
\small\tt zgao@math.carleton.ca\\}

\begin{document}
\maketitle

\begin{abstract}
We define the notion of asymptotically free for locally restricted
compositions, which means roughly that large parts can often be
replaced by any larger parts.
Two well-known examples are Carlitz and alternating compositions.
We show that large parts have asymptotically geometric distributions.
This leads to asymptotically independent Poisson variables for numbers of
various large parts.
Based on this we obtain asymptotic formulas for the probability of
being gap free and for the expected values of the largest
part, number of distinct parts and number of parts of multiplicity $k$,
all accurate to $o(1)$.
\end{abstract}

\centerline{\it Dedicated to the memory of Herb Wilf.}

\section{Introduction}

Various authors have considered aspects of unrestricted compositions
and Carlitz compositions (unequal adjacent parts) that require knowledge
about the large parts.
The results include information about largest part, number of distinct
parts, gap-freeness and number of parts of multiplicity $k$.
We extend these results to a broad class of compositions, drawing on
earlier work on {\em locally restricted compositions}~\cite{BC2} by
defining a subclass of locally restricted compositions for which
we can show that the large parts are asymptotically independent
geometric random variables.
This leads to asymptotically independent Poisson random variables
for numbers of various large parts.
Our main goal is to prove Theorem~\ref{thm:main}.
Although a full understanding of the theorem requires some definitions,
it can be read now.
Among the compositions included in our definition are unrestricted,
Carlitz and alternating up-down.

Although it was not possible to compute generating functions in~\cite{BC2},
various properties were established, including the following.
\begin{itemize}
\item[(a)]
The number of compositions of $n$ is
$Ar^{-n}(1+O(\delta^n))$ for some $0<\delta<1$ because of a simple
pole in the generating function. Since the convergence to
$Ar^{-n}$ is exponentially fast, the values of $r$ and $A$ can be estimated
fairly easily if one can count compositions for relatively small
values of $n$.
\cite[Theorem~3]{BC2}
\item[(b)]
If a subcomposition can occur
arbitrarily often, the number of times it occurs in a random
composition of $n$ has a distribution that is asymptotically
normal with mean and variance asymptotically proportional to $n$.
The same is true for the total number of parts in a random
composition.
\cite[Theorem~4]{BC2}
\item[(c)]
In many cases, the largest part and number
of distinct parts in a random composition is asymptotic to
$\log_{1/r}n$.
\cite[Section~9]{BC2}
\end{itemize}
Various special cases were considered ~\cite{BC1,BC3}, where more
could be said about the generating functions.
In none of these papers was the behavior of the large parts addressed beyond that
in~(c).

\begin{definition}
[Composition terminology]
\label{definition:composition terminology}
$\naturals$ and $\naturals_0$ denote the positive integers and the positive
integers and 0, respectively.

A composition is written $\cvec c=c_1\cdots c_k$ where $c_i\in\naturals$.
(We never write it as $c_1\ldots c_k$.)
We use the same notation to denote concatenation of compositions as
in $\cvec a_1\cdots\cvec a_m$.  The length of $\cvec c=c_1\cdots c_k$ is denoted by $\len(\cvec c)=k$
and the sum of the parts by $\partsum(\cvec c)$.

A {\em subcomposition} of $\cvec c$ is a sequence of one or more
consecutive parts of $\cvec c$.
The ordered $k$-tuple $(L_1,\dots,L_k)$ is a {\em subsequence of} $\cvec c$
if for some increasing sequence of
indices $1\le j_1<j_2<\cdots<j_k \le \len(\cvec c)$ we have
$c_{j_i} = L_i, 1\le i\le k$.
A subsequence of a composition is {\em marked} if the elements of
the subsequence are distinguished in some manner.
For example, in the composition {\it abacb} there is no marked subsequence
whereas {\it\.aba\.cb} and {\it ab\.a\.cb} each contain the
marked subsequence $(a,c)$.

\end{definition}

\begin{definition}
[Local restriction function]
\label{definition:LRF}
Let $m,p\in\naturals$.  A {\em local restriction function of type}
$(m,p)$ is a function
$$
\lrf:\{0,1,\dots,m-1\} \times (\naturals_0)^{p+1} \rightarrow \{0,1\}
$$
with $\lrf(i;0,\ldots,0)=1$ for all $i$.
The integers $m$ and $p$ are called, respectively, the {\em modulus}
and {\em span} of $\lrf$.
\end{definition}

\begin{definition}
[Class of compositions determined by local restrictions]
\label{definition:LRCPhi}
Let $\lrf$ be a local restriction function.
The {\em class of compositions determined by} $\lrf$ is
$$
\sC_{\lrf} = \{\cvec c: \cvec c {\rm ~is~a~composition,~and~}
                         \lrf(i\mod m;c_i,c_{i-1},\dots,c_{i-p})=1
                 {\rm ~for~} i \in \integers
             \}.
$$
If an index $j$ refers to a part before the first part ($j<1$)
or after the last part ($j$ greater than the number of parts),
we set $c_j=0$.

A class $\sC$ of compositions is
{\em locally restricted} if $ \sC=\sC_{\lrf} $
for some local restriction function~$\lrf$.

If $\lrfs$ is a set of local restriction functions, we define
$\sC_\lrfs = \cup_{\lrf\in\lrfs}\sC_\lrf$.
\end{definition}

\noindent
The number $m$ determines a periodicity.
The number $p$ determines a ``window''---by looking at parts $c_j$
with $0<|i-j|\le p$, we can determine what values, if any, are allowed
for $c_i$.
We could replace $m$ by any multiple of itself, $p$ by any
larger value, and redefine $\lrf$ to get the same class of compositions.

\BeginExample[Alternating compositions]\label{ex:up-down} Up-down
compositions $c_1\le c_2\ge c_3\le\cdots$ can be described as
follows. Set $m=2$, $p=1$,
\begin{equation}\label{eq:weak up-down}
\lrf(1;a,0)=1,\quad \lrf(1;a,b)=1,\quad
\lrf(0;0,a)=1\quad\mbox{and}\quad \lrf(0;b,a)=1,
\end{equation}
whenever $0<a\le b$. Otherwise set $\lrf=0$, except that
$\lrf(i;0,0)=1$ as required by Definition~\ref{definition:LRF}.

The function $\lrf$ describes alternating compositions that start
by going up (because $\lrf(1;a,b)=1$) and have an odd number of
parts (because $\lrf(0;0,a)=1$ permits the first zero after the
composition to be in an even position; but $\lrf(1;0,a)=0$ forbids
it to be in an odd position). We could have included an even
number of parts as well by defining $\lrf(1;0,a)=1$.

We cannot extend the definition of $\lrf$ to include compositions
that begin by going down. These can be defined by switching ~0
and~1 in the first argument of $\lrf$ to give a new function
$\lrf'$. With the extension to $\lrf$ (and hence $\lrf'$) noted in
the previous paragraph, $\sC_{\{\lrf,\lrf'\}}$ consists of all
alternating compositions.

Suppose we require that the inequalities be strict. This can be
done by simply changing $0<a\le b$ to $0<a<b$ in~(\ref{eq:weak
up-down}). Now, however, we can include all strict alternating
compositions in one $\lrf$ instead of using
$\sC_{\{\lrf,\lrf'\}}$. Set $m=1$, $p=2$,
$$
\lrf(0;0,0,a) = \lrf(0;a,0,0) = 1,\quad \lrf(0;0,a,b) =
\lrf(0;a,b,0) = 1\quad\mbox{and}\quad \lrf(0;a,b,c) = 1,
$$
when $a,b,c\in\naturals$ and either $a<b>c$ or $a>b<c$. It may
appear at first that $m=1$ causes periodicity to be lost; however,
by looking at the two previous parts we can determine which of
$c_{i-1}>c_{i-2}$ and $c_{i-1}<c_{i-2}$ holds. This will not work
with weakly alternating compositions since they can have
arbitrarily long strings of equal parts. \EndExample

\begin{definition}
[Recurrent compositions] Let $\sC$ be a class of locally
restricted compositions with span $p$ and modulus $m$.

We say that a subcomposition $\cvec s$ is {\em recurrent at $j$
modulo $m$} if, for every $k$ and every $\cvec a\,\cvec x\,\cvec
z\in\sC$ with $\len(\cvec a)\ge p$ and $\len(\cvec z)\ge p$, there
is a composition $\cvec a\cdots\cvec z\in\sC$ containing at least
$k$ copies of $\cvec s$ starting at positions congruent to $j$
modulo~$m$.
\begin{itemize}
\item If $\cvec s$ is recurrent for some $j$, we say $\cvec s$ is
{\em recurrent}. \item If a recurrent subcomposition has length 1,
we call it a {\em recurrent part}. \item A class $\sC_\lrf$ (and
$\lrf$) is {\em recurrent} if every subcomposition $c_i\cdots c_j$
of $\cvec c\in\sC$ with $i>p$ and $j+p\le\len(\cvec c)$ is
recurrent. \item A class $\sC_\lrfs$ (and $\lrfs$) is {\em
recurrent} if $\lrf$ is recurrent for every $\lrf\in\lrfs$.
\end{itemize}
\end{definition}

\noindent It is a consequence of these definitions that if $\cvec
r$ and $\cvec s$ are recurrent subcompositions, $\len(\cvec a)\ge
p$, $\len(\cvec z)\ge p$, and $\cvec a\,\cvec x\,\cvec z\in\sC$,
then there is a composition $\cvec a\cdots\cvec r\cdots\cvec
s\cdots\cvec z$ in $\sC_\lrf$. (We get $\cvec a\cdots\cvec r\cvec
y\cvec z$ for some $\cvec y$. Replace $\cvec a$ with $\cvec
a\cdots\cvec r$ and $\cvec x$ with $\cvec y$ in the definition.)

For the first $\lrf$ in Example~\ref{ex:up-down}, the 2-part
subcomposition $ab$ is recurrent at 1~modulo~2 whenever $0<a<b$
and is recurrent at 0~modulo~2 whenever $0<b<a$. The part 1 is
recurrent at 1~modulo~2 but is not recurrent at 0~modulo~2.

\BeginRemark[Ignoring nonrecurrent parts]
Since nonrecurrent parts can only appear in the first or last
$p$ parts, and since almost all compositions of $n$ have
$\Theta(n)$ parts, we can usually ignore
the nonrecurrent parts in our asymptotic estimates.
\EndRemark

\begin{definition}
[Similar restrictions]\label{def:similar}
Suppose $\lrf$ and $\lrf'$ are local restriction functions with the
same modulus and span.
Suppose $\sC_\lrf$ and $\sC_{\lrf'}$ are recurrent and there is a $k$
such that $\cvec s$ is recurrent at $j$~mod~$m$ in $\sC_\lrf$ if and
only if it is recurrent at $(j+k)$~mod~$m$ in $\sC_\lrf'$.
We then say that $\sC_\lrf$ and $\sC_{\lrf'}$ are {\em similar}
and write $\sC_\lrf\approx\sC_{\lrf'}$ as well as $\lrf\approx\lrf'$.
\end{definition}

\noindent
Clearly $\approx$ is an equivalence relation.

\BeginExample[Alternating compositions again]
In Example~\ref{ex:up-down}, $\lrf\approx\lrf'$ and so
it turns out that Theorem~\ref{thm:main} will apply to all
alternating compositions.
This remains true if we make either one or both of the
inequalities $x<y$ and $x>y$ weak.
However, weak and strong inequalities give restrictions
which are not similar.
For example, if the restrictions in $\lrf'$ were changed to
weak giving $\lrf''$, we would not have $\lrf\approx\lrf''$
and so we could not apply Theorem~\ref{thm:main} to
$\{\lrf,\lrf''\}$.
\EndExample

\BeginRemark[Some asymptotics]
We refer to (a) and (b) near the start of this section.
Since the radius of convergence $r$ in (a) depends only on
the recurrent subcompositions, it will follow that the form
$A(r^{-n}(1+O(\delta^n))$ still holds for $\sC_\lrfs$
when $\lrfs$ is a {\em finite} set of similar restrictions.
For essentially the same reason, the normality in (b) continues
to hold.
(See Section~\ref{sec:LRC2} for details.)
\EndRemark

\begin{definition}
[Asymptotically free]
\label{def:asymp free}
Let $\sC_\lrf$ be a set of locally restricted compositions
of span $p$.
If $\sC_\lrf$ is recurrent and
the following hold, we say that $\lrf$ and the compositions in
$\sC_\lrf$ are {\em asymptotically free.}
\begin{itemize}
\item[(a)]
Suppose $j$ and $r_i$ are such that
$\cvec r(x)=r_1\cdots r_pxr_{p+2}\cdots r_{2p+1}$
is recurrent at $j$ modulo $m$ for infinitely
many values of $x$.
Then there is an $M$ (depending on $j$ and the $r_i$)
such that, if $\cvec r(x)$ occurs at a position $j$
mod~$m$ in a composition, we may replace that $x$
by any $x'\ge M$.
\item[(b)]
There is at least one set of values $j$ and $r_i$ of the sort
described in (a).
\end{itemize}
Let $\lrfs$ be a {\em finite set of similar} local restriction functions.
If $\lrf\approx\lrf'$ and $\lrf$ is asymptotically free, then
clearly $\lrf'$ is asymptotically free.
Hence we say that $\lrfs$ and the compositions in $\sC_\lrfs$ are
{\em asymptotically free} if $\sC_\lrf$
is asymptotically free for some $\lrf\in\lrfs$.
\end{definition}

\noindent
It is fairly easy to verify that asymptotically free
$\sC_\lrf$ are special cases of the regular $\sC_\lrf$ studied in~\cite{BC2}.
Note that, since $\lrf$ has span $p$, no parts other than the $r_i$ impose
restrictions on $x$.
We arrived at the notion of asymptotically free as a concept succinctly
stated, fairly intuitive, and inclusive of a number of known examples, for
which the results of Theorem~\ref{thm:main} hold.
It would be of interest to extend these results to more classes of compositions.

\BeginExample[A bad definition]
We could have attempted to define asymptotically free $\lrfs$ by simply insisting
that~(a) and~(b) hold for $\sC_\lrfs$, however this is insufficient.
Consider $\lrfs=\{\lrf,\lrf'\}$ and $\lrf$ (resp.\ $\lrf'$) requires that parts in
odd (resp.\ even) positions be odd.
Then large odd parts will tend to be more common than large even parts and so the
conclusion in Theorem~\ref{thm:main}(a) would be false.
\EndExample

\BeginExample[Generalized Carlitz compositions]
Carlitz compositions are defined by the restriction $c_i\ne c_{i-1}$.
They were generalized to restricted differences in~\cite{BC1} by requiring that
$c_i-c_{i-1}\notin{\cal N}$ where $\cal N$ is an arbitrary set of integers.
(Carlitz compositions correspond to ${\cal N}=\{0\}$.)
These compositions are recurrent with modulus~1 and span~1.
If $\cal N$ is finite, we have asymptotically free compositions.
For the generalized Carlitz compositions studied in \cite{BC1},
$\cal N$  was the same for all $c_{i-1}$. We can generalize further
by letting  $\cal N$ depend on the value of $c_{i-1}$, say
${\cal N}(c_{i-1})$.  If all the ${\cal N}(c)$ are finite, we still have
asymptotically free compositions; however, they cannot be studied
by the method in~\cite{BC1}. Instead,~\cite{BC2} must be used.
\EndExample

\BeginExample[Some periodic conditions]
Up-down compositions have constraints of modulus~2.
General periodic inequality constraints were studied in~\cite{BC3}.
These are all asymptotically free provided they allow parts to
both increase and decrease.
As in the preceding example, we could require that the change between
adjacent parts be dependent on the parts. For example, we could require
that the ratio of adjacent parts be at least~2 ($c_i/c_{i-1}\ge2$ for an
increase and $c_{i-1}/c_i\ge2$ for a decrease).

For fixed $k$, $k$-rowed compositions $a_{i,j}$ in which
differences of adjacent parts avoid a finite set are
asymptotically free.
One interleaves the parts to produce a one-rowed composition:
If $a_{i,j}$ are the parts of a $k$-rowed composition of $n$, then
$c_{i+k(j-1)}=a_{i,j}$ for $1\le i\le k$ and
$j=1,2,\cdots$ gives a bijection with one-rowed compositions
$\cvec c$ of $n$.
We can take the modulus and span to be $k$.
\EndExample

\begin{definition}
[Gap free]\label{def:gap-free}
A composition with largest part $M$ is called {\em gap free} if it
contains all recurrent parts less than $M$.
\end{definition}

\noindent
The restriction of gap-free to recurrent parts is used to rule out
classes such as the following. Let $\sC$ be all compositions subject
to the restriction that 2 and 3 can appear only as the first
part of a composition. Since almost all compositions contain 1 and
no composition in $\sC$ can contain both 2 and 3, almost no
compositions in $\sC$ would be gap-free if we required that the
support of the parts be an interval in $\naturals$.

\begin{conventions}
We use the following conventions in this paper.
\begin{itemize}
\item
When we talk about something random, we always mean that it
is chosen uniformly at random from the set in question. We say that a property
holds {\em asymptotically almost surely} (a.a.s) if the probability that the property holds
tends to 1 as the size of the set goes to infinity, and we also say that
{\em almost all} objects in the set have the property.
\item
Expectation is denoted by $\Ex$.
\item
After a class of compositions has been defined, we usually omit
the modifiers (e.g.\  asymptotically free) and refer to elements of the
class simply as compositions.
\item
The number of compositions of $n$ in the class $\sC$ is
asymptotically $Ar^{-n}$.
We will always use $A$ and $r$ for these parameters.
\item
All logarithms are to the base $1/r$ except the natural logarithm $\ln$.
\end{itemize}
\end{conventions}

Remember that we call $\sC_\lrfs$ asymptotically free if and only if $\lrfs$
is a finite set of similar asymptotically free local restriction functions.
\begin{theorem}[Main theorem]\label{thm:main}
Let $\gamma\doteq 0.577216$ be Euler's constant and let
\begin{equation}\label{eq:oscillation}
P_k(x) ~= ~\log e\sum_{\ell\ne 0}
\Gamma(k+2i\pi \ell\log e)\exp(-2i\ell\pi \log x).
\end{equation}
(This is a periodic function of $\log x$.
For $1/2<r<1$ and $k=0$ the amplitude is less than $10^{-6}$.)

Let $\lrfs$ be asymptotically free and let $r$ be the radius of
convergence of the generating function for $\sC_\lrfs$.
The following are true for some $C>0$, which has the same value
in all parts of the theorem.
\begin{itemize}
\item[(a)]
Select a composition of $n$ uniformly at random.
Let $X_0(n)$ be the number of parts and $X_k(n)$ the number of parts
of size $k$.
For recurrent $k$ and $\epsilon>0$,
\begin{equation}\label{eq:Xk/X0}
{\rm Prob}\left(
\left|
   \frac{X_k(n)}{X_0(n)} - \frac{\Ex(X_k(n))}{\Ex(X_0(n))}
\right|
> \epsilon
\right)
 \to 0
~~\mbox{as $n\to\infty$}.
\end{equation}
Furthermore, the limit
\begin{equation}\label{eq:Ex/Ex}
\p_k = \lim_{n\rightarrow\infty}
\frac{\Ex(X_k(n))}{\Ex(X_0(n))}
\end{equation}
exists, and $\p_k \sim Br^k$ as $k\to\infty$ for some positive
constant $B$.
\item[(b)]
Let the random variable $M_n$ be the size
of the maximum part in a random composition of $n$. For any
function $\omega_b(n)$ such that $\omega_b(n)\to\infty$ as
$n\to\infty$, $|M_n-\log n|<\omega_b(n)$ a.a.s. Furthermore
$$
\Ex(M_n) ~=~
\log\left(\frac{Cn}{1-r}\right)+\gamma\log e
- \frac{1}{2}+P_0\!\left(\frac{Cn}{1-r}\right)+o(1),
$$
where $C = B\lim_{n\to\infty}\Ex(X_0(n))/n$.
\item[(c)]
Let $\nu$ be the number of nonrecurrent parts.
(Since the compositions are asymptotically free, $\nu$ is
finite.) Let the random variable
$D_n$ be the number of distinct $\underline{\mbox{recurrent}}$
parts in a random composition of $n$. For any function
$\omega_c(n)$ such that $\omega_c(n)\to\infty$ as $n\to\infty$,
$|D_n-\log n|<\omega_c(n)$ a.a.s. Furthermore
$$
\Ex(D_n) + \nu ~=~
\log(Cn)+\gamma\log e-\frac{1}{2}+P_0(Cn)+o(1).
$$
\item[(d)]
Let $q_n(k)$ be the fraction of compositions of $n$ which are
gap-free and have largest part $k$.
There is a function $\omega_d(n)\to\infty$ as $n\to\infty$ such that
\begin{equation}\label{eq:qn(k)}
q_n(k) ~\sim~
\exp\left(\frac{-Cnr^{k+1}}{1-r}\right)
\prod_{j\le k}\left(1-\exp\left(-Cnr^j\right)\right)
\end{equation}
uniformly for $|k-\log n|<\omega_d(n)$. Furthermore, for any
constant $D$, the minimum of $q_n(k)$ over $|k-\log n|<D$ is
bounded away from zero.
\item[(e)]
Let $q_n$ be the fraction of
compositions of $n$ which are gap-free. Then $q_n$ is asymptotic
to the sum of the right side of (\ref{eq:qn(k)}), where the sum
may be restricted to $|k-\log n|<\omega_d(n)$ for any
$\omega_d(n)\to\infty$ as $n\to\infty$. Furthermore, $q_n\sim p_m$
where $m=\left\lfloor\frac{Cn}{1-r}\right\rfloor$ and
\begin{equation}\label{eq:pm}
p_m ~=~ \cases{
1_{\vphantom{\bigm|}}& if~ $m=0$;\cr
\displaystyle\sum_{k=0}^{m-1}p_k{m\choose k}r^k(1-r)^{m-k}& if~ $m>0$.\cr}
\end{equation}
\item[(f)]
Let $g_n(k)$ be the fraction of compositions of $n$ that have exactly
$k$ parts of maximum size.
Then for each fixed $k$ and as $n\to\infty$,
$$
g_n(k) ~\sim~
\frac{(1-r)^k}{k!}P_k\left(\frac{Cn}{1-r}\right)+\frac{(1-r)^k\log
e}{k}.
$$
\item[(g)]
Let $D_n(k)$ be the number of distinct recurrent parts that appear
exactly $k$ times in a random composition of $n$. For fixed $k>0$
$$
\Ex(D_n(k)) ~=~ \frac{P_k(Cn)}{k!} + \frac{\log e}{k} + o(1).
$$
Let $m_n(k)$ be the probability that a randomly chosen recurrent part size
in a random composition of $n$ has multiplicity $k$.
For fixed $k$, $m_n(k)\sim \Ex(D_n(k))/\log n$.
\item[(h)]
Let $\lrfs'$ be a finite set of local restriction functions similar to
those in $\lrfs$.
The values of $r$, $B$ and $C$ are the same for $\sC_\lrfs$ and
$\sC_{\lrfs'}$.
\end{itemize}
\end{theorem}

\noindent
We recall that $\Gamma(a+iy)$ goes to zero exponentially fast as $y\to\pm\infty$.
Thus the sum (\ref{eq:oscillation}) is dominated by the terms with small $\ell$.

\smallskip

\noindent
Since estimating $C$ is generally harder than estimating $A$,
the following theorem is sometimes useful.

\begin{theorem}[Sometimes $A=C$]\label{thm:AeqC}
Let $\sC$ be a class of asymptotically free compositions
and let the number of compositions of $n$ be asymptotic
to $Ar^{-n}$.
Suppose that there is some $\ell$ such that,
whenever the number of parts in each of $\cvec a$ and
$\cvec b$ is at least $\ell$, we have that
$\cvec c=\cvec a x\cvec b$ is in $\sC$ for
infinitely many $x$ if and only if $\cvec a$ and $\cvec b$
are in $\sC$.
Then $C=A$, where $C$ is the constant in Theorem~\ref{thm:main}.
\end{theorem}

\begin{theorem}[Asymptotically Poisson]\label{thm:Poisson}
Let $\zeta_j$ be the number of parts of size $j$ in a random
composition in $\sC$ of size $n$.  Then there is a function
$\omega(n)\rightarrow\infty$ such that the random variables
$\{\zeta_j:\log n -\omega(n) \le j\le n\}$ are asymptotically
independent Poisson random variables with means $\mu_j =Cnr^j$.
\end{theorem}

\section{Discussion and Examples}

\BeginRemark
[Some previous results]
We review some results that involve the study of parts of large size.

Most results deal with unrestricted compositions.
As far as we know, the first result is due to Odlyzko and Richmond~\cite{OR}.
For $a(n,m)$, the number of compositions of $n$ with largest part $m$,
they prove the sequence is unimodal for each $n$ and show that the $m$
which maximizes $a(n,m)$ is always one of the two integers closest to $\log_2 n$.
The fact that the largest part $M_n$ is strongly concentrated is well known.
For example, it appears as an exercise in~\cite{FS}.
Hwang and Yeh~\cite{HY} studied the distinct parts in a random composition,
obtaining asymptotics for the expected value of their number and sum as well
as other results.
Hitczenko and Stengle~\cite{HS} also studied the expected number of distinct parts.
The asymptotic probability that a composition is gap-free was
obtained by Hitczenko and Knopfmacher~\cite{HK}.
They based their proof on a gap free result they obtained for samples of iid
geometric random variables, which we also use in our study of gap-freeness.
Wilf asked about $m_n(k)$, the probability that a randomly chosen part size in a
random composition of $n$ had multiplicity $k$.
This problem was studied by Hitczenko, Rousseau and Savage~\cite{HSv,HRS}.
Louchard~\cite{Lo} studied $D_n(k)$, obtaining information about its moments.
Archibald and Knopfmacher~\cite{AK} studied the largest missing part in
compositions that are not gap free.

Fewer results have been obtained for Carlitz compositions.
Using~\cite{FGD}, Knopfmacher and Prodinger~\cite{KP} obtained asymptotics
for the largest part in Carlitz compositions and observed that there was
oscillatory behavior.
The expected number of distinct parts, $\Ex(D_n)$, was studied by
Hitczenko and Louchard~\cite{HL} who required an independence assumption
that was eliminated by Goh and Hitczenko~\cite{GH}.
Kheyfets~\cite{Khe} obtains results for parts of multiplicity $k$ that
parallel those mentioned in the previous paragraph for $D_n(k)$ and
$m_n(k)$ in the unrestricted case.
Louchard and Prodinger~\cite{LP} study the distribution of part sizes.

Theorem~\ref{thm:main} extends most of these results to asymptotically free compositions.
One exception is~\cite{AK} which came to our attention when this paper was
essentially complete.
It is likely that our methods can generalize their results, although with less
accuracy than they obtain.
Most of the known results for unrestricted and Carlitz compositions have
greater accuracy than our results which typically have $o(1)$ error rather
than more explicit estimates.
Also, we do not have formulas for the two constants $C$ and $r$ appearing in
our results, whereas they are known for unrestricted and Carlitz compositions.
However, since the number of compositions is $Ar^{-n}$ with an exponentially
small relative error the more important $r$ is easily estimated if one can
count compositions for moderate values of $n$ efficiently.

An earlier version of this paper appeared, without proofs, as the extended
abstract~\cite{BCG}.
The present paper considers a more general class of compositions and
contains some additional results.
\EndRemark

\BeginExample[$A=C$]
It is easily seen that Theorem~\ref{thm:AeqC} applies to the following
classes of compositions
\begin{itemize}
\item[(a)]
unrestricted compositions (so $C=1/2$);
\item[(b)]
compositions where the value of $c_i$ is restricted only by
$c_{i-1}$ and $c_{i+1}$ and may be arbitrarily large;
\item[(c)]
alternating compositions ($c_{2i-1}<c_{2i}>c_{2i+1}$)
where the number of parts must be odd.
\end{itemize}
We note that (b) includes Carlitz compositions and so
$C=0.4563634741\cdots$ for Carlitz compositions~\cite{LP}.
The inequality conditions in (c) can be generalized:
we may require that $c_{2i}-c_{2i-1}$ and $c_{2i}-c_{2i+1}$ belong
to some subset of $\integers$ that contains arbitrarily large
positive values and the subset may depend on $i$ modulo some period.

Although (c) gives $A=C$ for only one type of alternating compositions,
it follows from Theorem~\ref{thm:main}(h) that the value obtained for~$r$,
$B$ and~$C$ in this case are the same for the various types of alternating
compositions discussed in Example~\ref{ex:up-down} even though they
have differing values of~$A$.
\EndExample

\BeginExample[Gap-free]\label{ex:gap-free}
The numbers $p_m$ in (e) were studied by Hitczenko and Knopfmacher~\cite{HK}
who showed that they oscillated with the same period as~(\ref{eq:oscillation})
when $r>1/2$.
They showed that, for $r=1/2$, there is no oscillation.
Their Figure~7 shows that the amplitude of oscillation of $p_m$ is less
than $10^{-6}$.
Consequently, if $r$ is known, one can determine the asymptotic value of $p_m$
and hence $q_n$ to within $10^{-6}$.
The following are the values of $p_m$ for three families of
compositions, correct up to the sixth decimal place.
\begin{itemize}
\item
For Carlitz compositions, it is known $r\doteq .57134979$.
It follows from (\ref{eq:pm}) that $p_m\doteq 0.372000$ for $m\ge 25$.
\item
For strictly alternating compositions ($c_{2i-1}<c_{2i}>c_{2i+1}$),
$r\doteq 0.63628175$ by~\cite{BC3}.
It follows from (\ref{eq:pm}) that $p_m\doteq 0.252277$ for $m\ge 25$.
\item
For weakly alternating compositions ($c_{2i-1}\le c_{2i}\ge c_{2i+1}$),
$r\doteq .57614877$ by~\cite{BC3}.
It follows from (\ref{eq:pm}) that $p_m\doteq 0.363144$ for $m\ge 25$.
\end{itemize}

Here is an alternative definition of gap-free based on the literature:
A composition is gap-free if, whenever it contains two recurrent parts,
say $a$ and $b$, it contains all recurrent parts between $a$ and $b$.
This definition does not alter the conclusions of
Theorem~\ref{thm:main}(d,e) because, by Lemma~\ref{lemma:normality}(b)
below, the fraction of compositions of $n$ that omit the smallest recurrent
part is exponentially small.
\EndExample

\BeginExample[Conjectures of Jakli\v c, Vitrih and \v Zagar]
Let ${\rm Max}_k(n)$ (resp.\ ${\rm Min}_k(n)$)
denote the number of all compositions of $n$ such that there are
more than $k$ copies of the maximal (resp.\ minimal) part.
Jakli\v c et al.~\cite{JVZ} conjectured that, when $k=1$
\begin{eqnarray}
& & \lim_{n\to \infty} \frac{{\rm Min_k}(n+1)}{{\rm Min_k}(n)} ~=~ 2
\label{eq:ConjMin}\\
& & \lim_{n\to \infty} \frac{{\rm Max_k}(n+1)}{{\rm Max_k}(n)} ~=~ 2.
\label{eq:ConjMax}
\end{eqnarray}
In fact, the conjectures hold for the compositions studied in this paper and
all $k\ge1$  provided 2 is replaced with $1/r$ and Min is restricted to recurrent parts.
The number of occurrences of any given recurrent part is $\Theta(n)$
for almost all recurrent locally restricted compositions of $n$ by~\cite{BC2}.
Thus (\ref{eq:ConjMin}) follows immediately from the fact that the number of
compositions of $n$ is asymptotic to $Ar^{-n}$.
We now prove (\ref{eq:ConjMax}).
Note that
$$
{\rm Max}_k(n) ~\sim~ Ar^{-n}\Biggl(1-\sum_{i\le k}g_n(i)\Biggr).
$$
By Theorem~\ref{thm:main}(f), $g_n(i)\sim g_{n+1}(i)$ and
$g_n(i)$ is bounded away from zero as $n\to\infty$.
Thus
$$
\frac{{\rm Max_k}(n+1)}{{\rm Max_k}(n)}
~\sim~
\frac{\Bigl(1-\sum_{i\le k}g_{n+1}(i)\Bigr)Ar^{-n-1}}
{\Bigl(1-\sum_{i\le k}g_n(i)\Bigr)Ar^{-n}}
~\sim~
\frac{1-\sum_{i\le k}g_{n+1}(i)}{1-\sum_{i\le k}g_n(i)}
\frac{1}{r}
~\sim~ \frac{1}{r}.
$$
One can change the definition of ${\rm Max}_k$ to mean exactly
$k$ copies of the maximal part and a similar proof will hold.
\EndExample

\BeginExample[Counterexamples without freeness]
It was shown in Theorem~1(f) of~\cite{BC1} that when differences of adjacent
parts are restricted to a finite set, the largest part is asymptotically almost surely of
order $\sqrt{\log n}$, so the bound in Theorem~\ref{thm:main}(a) fails.
\EndExample

\section{Statement of Lemmas}

The following six lemmas are used in our proofs of Theorems~
\ref{thm:main}, \ref{thm:AeqC} and~\ref{thm:Poisson}.

\begin{lemma}[Normality and tails]\label{lemma:normality}
Let $\sC_\lrfs$ be a class of  asymptotically free compositions
and let $d$ be arbitrary.
Let ${\cal R}$ be a possibly infinite nonempty set of recurrent
subcompositions each of which contains at most $d$ parts.
Assume that if we alter $\lrfs$ to forbid the elements of ${\cal R}$,
the resulting class of compositions is still recurrent.
Let the random variable $X_n$ be either the number of
occurrences of elements of $\cal R$ in a random composition of $n$
or the number of parts in a random composition of $n$.
The following are true.
\begin{itemize}
\item[(a)]
The distribution of $X_n$ is asymptotically normal with mean and
variance asymptotically proportional to $n$.
\item[(b)]
There are constants $C_i>0$ depending on what $X_n$ counts such that
$$
\Pr(X_n\!<\!C_1n) ~<~ C_2(1+C_3)^{-n}
~~\mbox{for all $n$.}
$$
\item[(c)]
Let $\cvec s$ be a subcomposition.
There is a constant $B$ dependent only on $\sC$ such that the
probability that a random composition contains at least one copy
of $\cvec s$ is at most $Bnr^{\partsum(\cvec s)}$.
\end{itemize}
\end{lemma}

\begin{definition}[The function $\free$]
\label{def:free}
As in Definition~\ref{def:asymp free} let $\cvec r(x)=r_1\cdots r_pxr_{p+2}\cdots r_{2p+1}$
where $\cvec r=\cvec r(0)=r_1\cdots r_p0r_{p+2}\cdots r_{2p+1}$.  For an
asymptotically free class the set
$$
S(\cvec r) = \{x: \cvec r(x) {\rm ~is~recurrent~}\}
$$
is either finite or co-finite.  So, there is a smallest integer $q(\cvec r)$ such that
$$
{\rm either~} [q(\cvec r),\infty) \subseteq S(\cvec r)
{\rm ~~or~~}
              [q(\cvec r),\infty) \subseteq \overline{S(\cvec r)}.
$$
Define
$$
\free(P) = \max\{q(\cvec r): r_i \le P {\rm ~for~} 1\le i\le 2p+1\}.
$$
\end{definition}

\noindent
It follows from the definition that if $\max(\cvec r)\le P$
and $x\ge\free(P)$ and $\cvec r(x)$ is recurrent at $j$ modulo $m$,
then $\cvec r$ is asymptotically free at $j$ modulo $m$.

\begin{definition}[$P$-isolated]
Suppose $\cvec c= c_{i-p}\cdots c_i\cdots c_{i+p}$
is a recurrent subcomposition.
If no $c_j$, except possibly $c_i$, exceeds $P$,
we call $c_i$ {\em $P$-isolated}.
\end{definition}

\noindent
A consequence of the definitions is that, whenever $x\ge\free(P)$
is $P$-isolated we are free to replace $x$ by any part that is of
size $\free(P)$ or greater.

\begin{lemma}[Large part separation]\label{lemma:separated}
Let $\sC_\lrfs$ be a class of asymptotically free compositions.
Suppose $\delta>0$.
There is a $P=P(\delta)$ and $N=N(\delta)$
such that the following holds for every $m\ge\free(P)$.
Let ${\cal M}(n)$ be the set of compositions of $n$ in which a part
of size $m$ has been marked.
For all $n>N+m$ the subset of ${\cal M}(n)$ in which the marked part
is not $P$-isolated has size less than $\delta|{\cal M}(n)|$.
\end{lemma}

\noindent
The following lemma proves most of Theorem~\ref{thm:main}(a).

\begin{lemma}[Geometric probabilities]\label{lemma:geometric}
We use the notation of Theorem~\ref{thm:main}(a).
\begin{itemize}
\item[(a)]
Equations (\ref{eq:Xk/X0}) and (\ref{eq:Ex/Ex}) are true.
\item[(b)]  Recall that $\p_k$ is the limit (on $n$)
of the ratio $\Ex(X_k(n))/\Ex(X_0(n))$.
For all sufficiently large parts $k$ and $\ell$ depending on $\delta>0$,
we have
$$
\left|\frac{\p_k\; r^{-k}}{\p_\ell\; r^{-\ell}}- 1\right| ~<~ \delta.
$$
\item[(c)]
We have $\p_k \sim Br^k$ for some positive constant $B$.
\end{itemize}
\end{lemma}

\begin{lemma}[Marked compositions]\label{lemma:L count}
Fix $k$ and a class $\sC$ of asymptotically free compositions.
Let $A$ be such that the number of compositions of $n$ is asymptotic to
$Ar^{-n}$ and let $C$ be as in Theorem~\ref{thm:main}.
If $\L(n)=(L_1(n),\ldots,L_k(n))$ is a sequence of $k$-tuples of
integers with
$$
\max(L_i)=o(n) ~~\mbox{and}~~ \min(L_i)\to\infty
~\mbox{as}~n\to\infty,
$$
then the number of compositions $\cvec c$ of $n$ having
$\L=(L_1 \cdots L_k)$ as a marked subsequence is
\begin{equation}\label{eq:L count}
(A+o(1))\frac{(Cn)^k r^{s-n}}{k!}
~\mbox{as}~ n\to\infty,
~\mbox{where}~ s=L_1+\cdots+L_k.
\end{equation}
\end{lemma}

\begin{lemma}[Characterization of Poisson]\label{GW12}
Let $[m]_k:=m(m-1)\cdots (m-k+1)$
denote the falling factorial.
Suppose that  $\zeta_1, \ldots, \zeta_n =\zeta_1(n), \ldots,
\zeta_n(n)$ is a set  of non-negative integer variables on a
probability space $\Lambda_n$, $n=1,2, \ldots,$ and there is a
sequence of positive reals $\gamma(n)$ and
  constants $0<\alpha<1$ and $0<c<1$ such that
\begin{description}
\item{(i)} $\gamma(n) \to \infty$ and $n-\gamma(n) \to \infty$;
\item{(ii)}  for any fixed positive integers $\ell$, $m_1,\ldots, m_\ell$,
and sequences
\hfill\break
\noindent $k_1(n)<k_2(n)<\cdots<k_{\ell}(n)$ with
$|k_i(n)-\gamma(n)|=O(1)$, $1 \le i \le \ell$, we have
\begin{equation} \label{eq:moments}
\Ex\left([\zeta_{k_1(n)}]_{m_1}[\zeta_{k_2(n)}]_{m_2}\cdots [\zeta_{k_\ell(n)}]_{m_\ell}\right)
\sim \prod_{j=1}^\ell \alpha^{(k_j(n)-\gamma(n))m_j},
\end{equation}
\item{(iii)}  $\Pr(\zeta_{k(n)}>0)=O\left(c^{k(n)-\gamma(n)}\right)$ uniformly
for all $k(n)>\gamma(n)$.
\end{description}
Then there exists a function $\omega(n)\to\infty$ so that
for $k=\lfloor\gamma(n)-\omega(n)\rfloor$, the total variation distance between
the distribution of $(\zeta_k,\zeta_{k+1},\ldots,\zeta_n)$, and that of
$(Z_k,Z_{k+1},\ldots,Z_n)$ tends to 0, where the $Z_j = Z_j(n)$ are
independent Poisson random variables with $\Ex Z_j = \alpha^{j-\gamma(n)}$.
\end{lemma}

\BeginRemark
The preceding lemma, Lemma~\ref{GW12}, is applied to obtain the Poisson result for large parts
stated as Theorem~\ref{thm:Poisson}.  The latter, in turn,
is used with Mellin transforms to prove Theorem~\ref{thm:main}(b-d);
and, with a result of Hitczenko and Knopfmacher~\cite{HK} on sequences of
geometric i.i.d.\ random variables, to prove Theorem~\ref{thm:main}(e,f).
\EndRemark

\begin{lemma}[Plentitude of recurrent parts]\label{lemma:zeta sum}
Let $\zeta_j$ be the number of occurrences of $j$ in
a random composition of $n$, and
let $k>0$ be arbitrary and fixed.
If $\omega(n)\to\infty$, then
$$
\sum_{\textstyle{j<\log n-\omega(n) \atop j{\rm ~recurrent}}}
\hskip-12pt
\Pr(\zeta_j\!<\!k) ~=~ o(1).
$$
\end{lemma}

\section{The Transfer Matrix and Sets of Functions}\label{sec:LRC2}

Before embarking on the proofs, we summarize some facts from~\cite{BC2}
which will be used and reduce the study of a finite set $\lrfs$ to a
single $\lrf$ since only single $\lrf$'s were considered in~\cite{BC2}.

We may replace the span $p$ by any larger value without altering the set
of compositions, provided we adjust the definition of $\lrf$.
Thus {\em we will assume that the span is a multiple of the modulus $m$.}
(Refer back to Definition~\ref{definition:LRF} for terminology.)

Let $\NC(n)$ be  the number of compositions of $n$ in a regular,
locally restricted class $\sC_\lrf$, and let $F(x)=\sum\NC(n)x^n$ be the ogf
(ordinary generating function).
Then, as proven in Theorem 2 of
that paper,
\begin{equation}\label{eq:LR2gf}
F(x^2) = \varphi(x) + F_{NR}(x^2),
\end{equation}
where
\begin{equation}\label{eq:varphi}
\varphi(x)= {\bf s}(x)^{\Tr} \, \Biggl( \sum_{k=0}^{\infty}T(x)^k \Biggr) \, {\bf f}(x).
\end{equation}
We now explain the various parts of~(\ref{eq:LR2gf}).
Here a ``small number of parts'' is at most some small multiple of $p$.

\noindent{\bf The transfer matrix $T(x)$}
is defined in terms of a certain sequence of
words $\cvec\nu_1, \cvec\nu_2, \dots$, where by a {\it word}
we mean a recurrent subcomposition of length $p$, the span,
(see Definition~\ref{definition:LRF}) whose first part is at j~mod~$p$
where $j$ is the same for all words indexing $T$ (and thus ${\bf s}$
and ${\bf f}$ as well).
The list contains all such recurrent words and
$$
T(x)_{ij} ~=~ \cases{
x^{\Sigma(\cvec\nu_i)+\Sigma(\cvec\nu_j)} &
if $\cvec\nu_j$ can follow $\cvec\nu_i$,\cr
0 & otherwise.}
$$
Except for parts near the ends, every composition is a concatenation of such words.
A single application of the transfer matrix corresponds to the adjunction of $p$
additional parts to the composition.

\noindent{\bf The infinite vectors ${\bf s}(x),{\bf f}(x)$}
have analytic entries corresponding to compositions with a small number of parts.
The component $s_i(x)$ of ${\bf s}(x)$ deals with the generating function for
the beginning of compositions where the last $p$ parts in the beginning are $\cvec\nu_i$.
Similarly, $f_j(x)$ deals with the generating function for endings whose first
$p$ parts are $\cvec\nu_j$.

\noindent{\bf The function $F_{NR}(x)$}
is the ogf for the subclass of compositions not counted in $\varphi(x)$.
These compositions have at most some small number of parts.
The ogf $F_{NR}(x)$ has radius of convergence 1.
(This is slightly different from the definition of $F_{NR}$ in~\cite{BC2};
however, all that matters for the theory is that $F_{NR}$ has radius of
convergence~1 and that~(\ref{eq:LR2gf}) counts all compositions exactly once.)

To assure that $T(x)$ satisfies certain useful technical conditions, it is
necessary to have the arguments $x^2$ and $x$ as indicated in~(a).
See the latter part of this section and~\cite{BC2} for more details on $T(x)$.

\subsection{Reduction to a single $\lrf$ and Theorem~\ref{thm:main}(h)}

Before discussing the more technical issues related to asymptotics, we explain
why it suffices to consider one $\lrf$ instead of an entire finite set $\lrfs$ of
similar $\lrf$.
This discussion will also prove Theorem~\ref{thm:main}(h).

Suppose $\lrf\approx\lrf'$ and let $T$ be the transfer matrix for $\lrf$.
Since $T_{ij}\ne0$ if and only if $\nu_i\nu_j$ is recurrent, we can use the
same transfer matrix for $\lrf'$; however, the vectors ${\bf s}$ and ${\bf f}$
will be different.
In fact, if $k$ is as in Definition~\ref{def:similar}, the number of parts
in the subcompositions of the vectors ${\bf s}$ for $\lrf$ and $\lrf'$
will differ by $k$~mod~$m$.
Nearly all results in~\cite{BC2} depend on $T$ but not on ${\bf s}$ or ${\bf f}$.
The exception is the constant $A$ in the asymptotic estimate $Ar^{-n}$ for the
number of compositions of $n$.

It follows that, if the sets $\sC_\lrf$, $\lrf\in\lrfs$, were pairwise disjoint
we could simply obtain results for one $\lrf\in\lrfs$ and combine the results where,
whenever $A$ is present we simply sum the values of $A$ for the various
$\lrf\in\lrfs$.
We now show that this can, in principle, be done.
There is no need to do this in practice since analytic methods for obtaining reasonable
estimates of $A$ are seldom available even for a single $\lrf$.

Fix temporarily a $\lrf\in\lrfs$.
Partition the elements $\lrf'$ of $\lrfs$ into $m$ sets $\lrfs_0,\ldots,\lrfs_{m-1}$
according to the value of $k$ in Definition~\ref{def:similar}.
We now focus on these sets,
first considering functions in different sets and then functions in the same set.

Suppose $\lrf\in\lrfs_i$ and $\lrf'\in\lrfs_j$ where $i\ne j$.
Consider the compositions in $\sC_\lrf\cap\sC_{\lrf'}$.
Let the value of $\lrf\lrf'$ be simply the product of $\lrf$ and $\lrf'$.
We note that $\sC_\lrf\cap\sC_{\lrf'}=\sC_{\lrf\lrf'}$ since a composition is in the
intersection if and only if it satisfies both local restriction functions.
Since $i\ne j$, it follows that the transfer matrix for the intersection will be
the same as that for $\lrf$ with some nonzero entries replaced by zeroes.
By Lemma~2(f) of~\cite{BC2} and the realization that the spectral radius
determines the growth rate (see below) it follows that the number of compositions
of $n$ in the intersection grows at an exponentially smaller rate than the number
in $\sC_\lrf$.
Hence, for asymptotic purposes, we may treat the $m$ sets $\sC_{\lrfs_i}$ as if
they are disjoint.
It follows that, except for Theorem~\ref{thm:AeqC}, we may assume we are dealing
with just one $\lrfs_i$.

We now consider a single $\lrfs_k$.
Suppose $\cvec c\in\sC_{\lrfs_k}$ is counted by (\ref{eq:varphi}).
We can write it in the form $\cvec a\cvec b\cvec z$, where $\cvec b$
is a sequence of words $\cvec\nu$ that index $T$, $\cvec s$ and $\cvec f$,
$\cvec a$ ends with one of these $\nu$ and $\cvec z$ starts with one of them.
By absorbing a recursive word or two in $\cvec a$ and $\cvec z$ if needed,
we can insure that the following two assumptions hold for some $\ell_i$.
\begin{itemize}
\item[(i)]
Since all compositions come from the same $\sC_{\lrfs_k}$ we can assume
$\len(\cvec a)=\ell_0$, the same value for all compositions in $\sC_{\lrfs_k}$.
\item[(ii)]
Since multiplication by $T$ adds $p$ parts to the compositions,
we can assume that the longest and shortest values of $\len(\cvec z)$,
say $\ell_1$ and $\ell_2$, differ by less than $p$.
\end{itemize}
It follows that each composition in $\sC_{\lrfs_k}$ with at least $\ell_0+\ell_1$
parts is counted by ${\bf s}T^k{\bf f}$ for some $k$ and has uniquely determined
$\cvec a$ and $\cvec z$.
We can limit attention to compositions with at least $\ell_0+\ell_1$ parts since
the generating function for those with fewer parts has radius of convergence at
least~1.

With each $\lrf_i\in\lrfs_k$ we associate two sets $\set S_i$ and $\set F_i$ as
follows.
$\cvec a\in\set S_i$ and $\cvec z\in\set F_i$ if and only if they satisfy
(i) and (ii) above and $\cvec a\cvec b\cvec z\in\sC_{\lrf_i}$ for some $\cvec z$.
The set $\set S_i$ determines $\bf s$ as follows.
If $\cvec a\in\set S_i$ ends with $\nu_j$, then a generating function obtained
from $\cvec a$ is added to $s_j$.
A similar construction holds for $\cvec z$ and $\bf f$.
Thus $\set S_i\times\set F_i$ determines the compositions in $\sC_{\lrf_i}$.
If we had $(\set S_i\times\set F_i)\cap(\set S_j\times\set F_j)=\emptyset$,
it would follow that $\sC_{\lrf_i}\cap\sC_{\lrf_j}$ would contain at most
some compositions shorter than $\ell_0+\ell_1$.
Thus, we need only prove that a union of Cartesian products
$\cup_{\lrf_i\in\lrfs_k}(\set S_i\times\set F_i)$
can always be written as a disjoint union of such products.  This is done by
considering the given terms $\set S_i\times\set F_i$ one at the time, and
using the identity
$$
(A\times B) \cap (C \times D)^{c} ~=~ \Bigl(A\times(B\setminus D)\Bigr)
\cup \Bigl((A\setminus C)\times (B\cap D)\Bigr),
$$
where the union is disjoint.
(Think of $C\times D$ as the latest $\set S_i\times\set F_i$, and
$A\times B$ as one of the pairwise disjoint components of the
previously processed $(i-1)$ products.  We keep $C\times D$ as a new
component, and each previously existing component is replaced by two
disjoint pieces.)
For each product in the resulting disjoint union, we construct a $\lrf$,
and their sum is the generating function for $\sC_k$, with the
possible exception of short compositions.




\subsection{Analytic aspects of $T(x)$ from~\cite{BC2}}

By Lemma 3 of~\cite{BC2}, at each $x_0\in(0,1)$
we have a neighborhood and functions $\lambda(x), E(x), B(x)$
analytic in that neighborhood such that
\begin{equation}
\label{eq:decomposition}
T(x)=\lambda(x) E(x)+B(x),
\end{equation}
where $E(x)$ is the projection onto the one-dimensional eigenspace
of eigenvalue $\lambda(x)$, the spectral radius of $B(x)$ is less than $\lambda(x)$,
and $E(x)B(x)=B(x)E(x)=0$.  The proof of Lemma 3 relies heavily on results
and methods from \cite{Kato}.
If we choose for $x_0$ the point $r^{1/2}$, $0<r<1$, where $\lambda(r^{1/2})=1$, it
follows from (\ref{eq:decomposition}) that
\begin{equation}\label{eq:LR2import}
\sum_{i\ge 0} T(x)^i ~=~ \frac{\lambda(x)}{1-\lambda(x)}E(x) +
(I-B(x))^{-1},
\end{equation}
in a punctured neighborhood $0<|x-r^{1/2}|<\delta$.  (We use $r^{1/2}$ so that $r$
is the radius of convergence for the ogf $F(x)$, and we are consistent with
our convention that $C(n)\sim Ar^{-n}$.)

The neighborhood $|x-r^{1/2}|<\delta$, which we shall refer to as
the $\delta$-neighborhood, plays a key role in our proof of
Lemma~\ref{lemma:L count}.  At any point in this neighborhood,
except the center $x=r^{1/2}$, the relation (\ref{eq:LR2import})
holds. At any $x_0$ with $|x_0|=r^{1/2}$, except $x_0=\pm r^{1/2}$,
the spectral radius of $T(x_0)$ is strictly less than 1 by Lemma 1
of \cite{BC2}.  Near such an $x_0$, the sum $S(x)=\sum_{i\ge
0}T(x)^i$ converges and $S(x)$ is analytic in a neighborhood of
$x_0$.  On the other hand, near $x_0=r^{1/2}$ the equation
(\ref{eq:LR2import}) shows that still $S(x)$ is analytic except for
an isolated singularity at $x=r^{1/2}$.  It is also shown in
\cite{BC2} that the root of $\lambda(x)=1$ at $x=r^{1/2}$ is a
simple root; thus near $r^{1/2}$ we have an analytic $\beta(x)$ with
$\lambda(x)=1-\beta(x)(1-x/r^{1/2})$, and $\beta(r^{1/2})\neq 0$.

\section{Proof of Lemma~\ref{lemma:normality}}

Part (a) follows from \cite[Thm.~4]{BC2} with just one
random variable $Y_1(n)=X_n$.
(The definition of ``unrelated events'' for that theorem is somewhat technical.
The condition in Lemma~\ref{lemma:normality} that altered $\lrfs$ be
recurrent insures that it holds.)

We now prove (b).
Into the transfer matrix $T(x)$ of~\cite{BC2}, introduce a new
variable $0<s\le1$ that keeps track of the number of
occurrences of elements of $\cal R$ or simply the number of parts.
If $d$ exceeds the span of $\lrf$, $d$ behaves like a new span and
it will be necessary to change $T$ so that one application of $T$
adds more parts to the composition.
Call the new matrix $T(x,s)$ and call the largest eigenvalue
$\lambda(x,s)$.
This leads to the asymptotics $A(s)r(s)^{-n}$ where $r(s)$
is the solution to $\lambda(r^{1/2}(s),s)=1$.
The case $s=1$ corresponds to the asymptotics for $\NC(n)$,
the number of compositions of $n$.
In the general case, we have asymptotics for
$\sum_k \NC(n,k)s^k$ where $\NC(n,k)$ is the number of compositions
of $n$ with exactly $k$ copies of elements of $\cal R$.
It follows that
$$
\sum_{k<\delta n} \NC(n,k)
~\le~ s^{-\delta n}\sum_k \NC(n,k)s^k
~\sim~ A(s)(s^\delta r(s))^{-n}.
$$
Hence it suffices to show that
\begin{equation}\label{eq:lambda s}
s^\delta r(s) > r(1)
\end{equation}
for some $s$ and $\delta$.
By Lemma~2(f) of~\cite{BC2}, $\lambda(x,s)<\lambda(x,1)$ for $x>0$ and $0<s<1$.
Since $\lambda$ is monotonically increasing in $x$, $r(s)>r(1)$ and so
(\ref{eq:lambda s}) holds for all sufficiently small $\delta$ depending on $s$.
This completes the proof of (b).

The proof of (c) is essentially the same as that given in
Section~9 of~\cite{BC2} for large part size. Since there are
slight changes, we repeat it here for completeness.

Let $p$ be the span of $\lrf$. Consider an expanded class
$\sC_\Psi$ where $\Psi$ is the same as $\lrf$ except that the
first $p$ and last $p$ parts of compositions are unrestricted.
The transition matrix $T(x)$ is unchanged.
Therefore $\Psi$ has the same radius of convergence $r$
as $\lrf$. Hence the number of compositions of $n$ in $\sC_\Psi$
is bounded above by $Cr^{-n}$ for some $C$. Hence the generating
function for $\sC_\Psi$ compositions by sum of parts is bounded
coefficient-wise by $C(1-x/r)^{-1}$.

Imagine marking a copy of $\cvec s$ in each composition in
$\sC_\lrf$. By the previous paragraph, the generating function for
such compositions of $n$ is bounded coefficient-wise by
$$
\frac{C}{1-x/r}x^s\frac{C}{1-x/r}
\mbox{~~where $s=\partsum(\cvec s)$.}
$$
Hence the number of such compositions of $n$ is bounded above by
$nC^2r^{s-n}$.

The previous paragraph overcounts the number of compositions
containing $\cvec s$. For some $C'>0$, the total number of
compositions of $n$ is at least $C'r^{-n}$ for large $n$. Taking
the ratio gives (c) with $B=C^2/C'$.
\EndProof

\section{Proof of Lemma~\ref{lemma:separated}}

Throughout the proof, whenever a new
condition is imposed on $P$ or $N$ it is understood that the
implied values must be at least as large as those already chosen.
All implied limits, as in $o(1)$, are as $n\to\infty$.

Let ${\cal M}^*(n)$ be the subset of ${\cal M}(n)$ in which the
marked part is not $P$-isolated.
Let $M^*(n)$ and $M(n)$ be the cardinalities of these two sets.
We will overestimate $M^*(n)$ and underestimate $M(n)$ and
show that their ratio can be made arbitrarily small provided $P$
and $n-m$ are sufficiently large.

For both counts, we consider compositions of the form
$\cvec a\cvec b\cvec c$ where $\cvec b$ contains a
special sequence of parts.

We then sum over $a$.
The composition $\cvec a$ will be like compositions in the
class $\sC_\lrf$ except that there will
be conditions on the last $p$ parts.
Since $T(x)$ is unchanged, the radius of convergence is
unchanged and so the number of $\cvec a$ is $\Theta(r^{-a})$
as $a\to\infty$.
A similar result holds for $\cvec c$.
We refer to this below as ``theta''.
Let $a=\partsum(\cvec a)$ and $c=\partsum(\cvec c)$.

We start with the underestimate of $M(n)$.
Let $\cvec r(x)$ be as in Definition~\ref{def:asymp free}
and let $b$ be the sum of its parts {\em excluding} $x$.
Choose $P$ so that $x$ is $P$-isolated in $\cvec r(x)$.
Thus $b$ is fixed as $n\to\infty$,
but we may increase $P$ as necessary later.
Since $m\ge\free(P)$, we may replace $x$ with a marked part $m$.
Let $\cvec b=\cvec r(m)$.
To underestimate $M(n)$, we will obtain a lower bound
on the number of occurrences of $\cvec b$.
Since the $r_i$ are fixed and the span is $p$,
the choices for $\cvec a$ and $\cvec c$ such that
$\cvec a\cvec b\cvec c\in\sC$ are independent of $m$.
By theta there are $B$ and $s$ such that there are at least
$Br^{-a}$ choices for $\cvec a$ and $Br^{-c}$ for $\cvec c$ when
$a\ge s$ and $c\ge s$.
Thus the total number of occurrences of $\cvec b$ is at least
$$
\sum_{a=s}^{n-(b+m)-s}B^2r^{-s}r^{-(n-s-(b+m))} ~=~
(n-2s-b-P+1)B^2r^{-n+b+m},
$$
and so for sufficiently large $n$ and some constant $C_0<B^2r^b$,
$$
M(n) ~\ge~ nC_0r^{-n+m}.
$$

For the non-isolated overcount, let $\cvec a$ and $\cvec c$ be
compositions where we put no restrictions on how they begin or
end. By theta the number of such compositions of $\ell$ is bounded
above by $Dr^{-\ell}$ for some $D$.
The composition $\cvec b$ will
contain at most $p+1$ parts.
It will either begin or end with the marked part $m$ and the other
ending part will be at least $P$ so that the marked part is not
$P$-isolated.
Let $b$ be the sum of the parts in $\cvec b$, omitting the marked part $m$.
It follows that $b\ge P$.
We bound the number of $\cvec b$ as follows.
Ignore the part $m$.
Choose a first part $b_1$ in $b$ ways.
Choose an additional $p-1$ parts, allowing parts of size zero, which
will be ignored when constructing $\cvec b$.
Since the remaining parts sum  to $b-b_1\le b-1$, each of them has at
most $b$ values.
Hence we have the bound $2\cdot b\cdot b^{p-1}$, where the factor of 2
is arises from the choice of which end to place $m$.
Thus, for some constant $C_1$,
$$
M^*(n) ~\le~ C_1\sum_{b=P}^n\sum_{a=0}^{n-b-m} r^{-a}b^p r^{-(n-a-(m+b))}
~<~ C_1nr^{-n+m}\sum_{b\ge P}b^pr^b.
$$

Combining our two estimates, we have for some constant $C_2$
and sufficiently large $n$
$$
\frac{M^*(n)}{M(n)} ~<~ C_2\sum_{i\ge P}i^p r^i +o(1).
$$
By choosing $P$ sufficiently large, we can make this arbitrarily small.
\EndProof

\section{Proof of Lemma~\ref{lemma:geometric}}

By Lemma~\ref{lemma:normality}(a), the total number of parts and the
number of parts of size $k$ are asymptotically normally distributed with
means and variances proportional to $n$.
Thus (a) follows.

We now prove (b).
Let $p$ be the span of $\lrf$.
Apply Lemma~\ref{lemma:separated} to obtain $P=P(\delta')$,
where $\delta'$ is sufficiently small and depends on the value of
$\delta$ in~(b).

Choose $k$ and $\ell$ larger than $\free(P)$.
Later we let $P$ and hence $k$ and $\ell$ tend to infinity slowly.

Consider compositions of $n+k$ with a marked part of size $k$.
By changing a part of size $k$ into one of size $\ell$ we obtain a
composition of $n+\ell$ with a marked part of size $\ell$.
This is a bijection between compositions containing a marked $P$-isolated part
of size $k$ and those containing a marked $P$-isolated part of size $\ell$,
the marked part being the one that is changed.
By Lemma~\ref{lemma:separated}, we can ignore those compositions with marked
parts that are not $P$-isolated.
Since the number of compositions of $m$ is asymptotic to $Ar^{-m}$,
the number of compositions with such $P$-isolated marked parts is asymptotic
to both $\Ex(X_k(n+k))Ar^{-n-k}$ and $\Ex(X_{\ell}(n+\ell))Ar^{-n-\ell}$.
Since $\Ex(X_0(n+k))\sim\Ex(X_0(n))\sim\Ex(X_0(n+\ell))$ as $n\to\infty$
with $k=o(n)$ and $\ell=o(n)$, (b) follows.

Part (c) follows by letting $\ell\to\infty$ in (b):
Since $\delta$ can be made arbitrarily small by choosing $k$ sufficiently
large, it follows that $\lim \p_\ell r^{-\ell}$ must exist and be nonzero.
\EndProof

\section{Proof of Lemma~\ref{lemma:L count}}

Let $Q(n)=\free(P(n))$ where $P(n)$ is some unspecified value that we will allow to increase
``sufficiently slowly'' with $n$.
When referring to $P$, $Q$ and the $L_i$ in the statement of the lemma, we will omit ``$(n)$''.
$P$ must increase so slowly that $\min(L_i)\ge Q$.
Let $s=L_1+\cdots+L_k$ and $m=n-s+kQ$.
Since $L_i=o(n)$, we have $m\sim n$.

Let $\bfR=(R_1,\dots,R_k)$ denote an arbitrary $k$-tuple of positive integers.
$$\mbox{
Denote the $k$-tuple $\bfR$ with $R_i=Q$ for $1\le i\le k$ by $Q^k$.
       }
$$
Let $\calM(n,\bfR)$ be the set of compositions of $n$ that have $\bfR$ as a marked
subsequence and let $M(n,\bfR)=|\calM(n,\bfR)|$.
We would like to establish a bijection between $\calM(n,\L)$ and $\calM(m,Q^k)$
by simply replacing the elements of one marked subsequence with those of the other.
Unfortunately this may fail if any of the following hold:
\begin{itemize}
\item[(a)]
an element of the marked subsequence has a part exceeding $P$ within distance $p$;
\item[(b)]
an element of the marked subsequence occurs within the first $p$ parts;
\item[(c)]
an element of the marked subsequence occurs within the last $p$ parts.
\end{itemize}
(The reason for (b) and (c) is that $\free$ applies only to the recurrent parts
of the composition and the ends may not be recurrent.)
Let the subscript $*$ refer to those compositions for which none of
(a)--(c) hold, {\em except that $p$ is replaced by $2p$ in (c)}.
Note that the proposed bijection is actually a bijection when restricted to
$\calM_*(n,\L)$ and $\calM_*(m,Q^k)$.
We will show that
\begin{equation}\label{eq:M*M}
M_*(m,Q^k)\sim M(m,Q^k) ~~\mbox{and}~~ M_*(n,\L)\sim M(n,\L).
\end{equation}
It then follows that
\begin{equation}\label{eq:MvsM}
M(n,\L)\sim M_*(n,\L)=M_*(m,Q^k)\sim M(m,Q^k)
\end{equation}
and so it suffices to estimate the size of any set that contains
$\calM_*(m,Q^k)$ and is contained in $\calM(m,Q^k)$.

\medskip

\noindent{\bf Overcounting compositions in $\calM\setminus\calM_*$}.
The idea is to allow the parts within distance $p$ of a marked part to be arbitrary.
We separate the composition into a sequence of $k+1$ possibly empty
subcompositions by removing the $k$ marked parts.
Each of the subcompositions so obtained may have been shifted with
regard to its modulus and may have beginning and ending subsequences
that are not allowed by the local restriction $\lrf$.
We overcount them by counting compositions with arbitrary shifts
in their moduli and with the initial and final $p$ parts arbitrary.
The generating function for these possible subcompositions has the
form (\ref{eq:LR2gf}); however, the values of ${\bf s}(x)$ and
${\bf f}(x)$ will be different.
Since $T(x)$ is unchanged, it is a consequence of the results
in~\cite{BC2} that the number of such subcompositions of $i$
is bounded above by $K_1r^{-i}$ for some $K_1$ and so the generating
function is bounded coefficient-wise by $K_1(1-x/r)^{-1}$.

The generating function for compositions that have a nonisolated
marked part $Q$ is bounded coefficient-wise by
$$
\left(\frac{K_1}{1-x/r}\right)^{k+1}\!k\,(2px^P)x^{kQ},
$$
where
\begin{itemize}
\item
the first factor bounds the subcompositions,
\item
$k$ chooses a marked part,
\item
$2p$ bounds the choices of a part near the chosen marked part,
\item
$x^P$ increases that part by $P$ to insure that a nearby part exceeds $P$, and
\item
$x^{kQ}$ inserts the marked parts $Q^k$.
\end{itemize}
The coefficient of $x^m$ is bounded by $K_2m^kr^{-m+kQ+P}$ for some constant $K_2$.
This takes care of all cases except the occurrence of two nearby marked parts $Q$.

Suppose there are two nearby marked parts.
Modify the previous argument by removing these two nearby parts and all
parts between them.
The termwise bound is now given by the generating function
$$
\left(\frac{K_1}{1-x/r}\right)^k
(k-1)\Biggl(\sum_{j=0}^p\frac{x}{1-x}\Biggr) x^{kQ},
$$
where the $k-1$ chooses a position to insert the pair of nearby
marked $Q$'s and the summation inserts arbitrary parts between these
two $Q$'s. Since the summation has radius of convergence~1 and
$r<1$, the coefficient of $x^m$ is bounded by $K_3m^{k-1}r^{-m+kQ}$.
Thus
\begin{eqnarray}
M(m,Q^k) - M_*(m,Q^k)
&=& m^kr^{-m+kQ}\Bigl(O(r^P) + O(1/m)\Bigr)
\nonumber\\
&=& n^kr^{-n+s}\Bigl(O(r^P) + O(1/n)+o(1)\Bigr)
\label{eq:Mdiff}
\end{eqnarray}

When a composition in $\calM(n,\L)$ is transformed by replacing $\L$ by
$Q^k$ and the result is an illegal composition, it must be of the
form we have just bounded and so the bound in (\ref{eq:Mdiff}) is
also a bound for $M(n,\L)-M_*(n,\L)$.

\medskip

\noindent{\bf Building marked compositions}.
We now build marked compositions that form a set between $\calM_*(m,Q^k)$
and $\calM(m,Q^k)$.
Define the transfer matrix $A_Q(x)$ by
$$
A_Q(x)_{i,j}~=~ \cases{
T(x)_{i,j}& if $\nu_i$ has exactly one $Q$;\cr
0 & otherwise.}
$$
Let $S(x)=\sum_{i\ge0}T(x)^i$.
Define the power-series $f_{k,Q}(x)$ by
\begin{equation}
\label{eq:Mgf}
   f_{k,Q}(x^2) =
       {\bf s}(x)^{\Tr} \, \left(S(x) A_{Q}(x)\right)^k
       \, S(x) \, {\bf f}(x).
\end{equation}
Since $[x^n]f_{k,Q}(x)$ equals the number of compositions counted by
$M(n,Q^k)$ having at most one marked $Q$ per word and no marked $Q$'s
near the ends of the composition, it follows that
$$
M_*(n,Q^k) \le [x^n]f_{k,Q}(x) \le M(n,Q^k).
$$
Using (\ref{eq:LR2import}) in (\ref{eq:Mgf}):
\begin{eqnarray}
f_{k,Q}(x^2) &=&
{\bf s}(x)^\Tr
\left(\frac{\lambda(x)}{1-\lambda(x)}E(x)A_Q(x) + B_Q(x)\right)^k
\nonumber\\
&&\null\times
\left(\frac{\lambda(x)}{1-\lambda(x)}E(x) + (I-B(x))^{-1}\right){\bf f}(x),
\label{eq:Mproduct}
\end{eqnarray}
where $B_Q(x) = (I-B(x))^{-1}A_Q(x)$.
When the products in (\ref{eq:Mproduct}) are expanded, we obtain something of the form
$$
f_{k,Q}(x^2) ~=~
\left(\frac{\lambda(x)}{1-\lambda(x)}\right)^{k+1}
{\bf s}(x)^\Tr(E(x)A_Q(x))^kE(x){\bf f}(x) + \frac{h_Q(x)}{(1-\lambda(x))^k},
$$
where $h_Q(x)$ is analytic in a neighborhood of $r^{1/2}$ because everything
in~(\ref{eq:Mproduct}) except $\frac{\lambda(x)}{1-\lambda(x)}$ is.
The first term determines the leading asymptotic behavior of the coefficients.
Recalling that $E(x)$ is a projection onto the 1-dimensional eigenspace of $\lambda(x)$,
define the functions ${\bf v}(x)$ and $\alpha_Q(x)$ by the equations
$$
{\bf v}(x)=E(x){\bf f}(x), {\rm ~~and~~} E(x) A_Q(x) {\bf v}(x) = \alpha_Q(x) {\bf v}(x),
$$
which are analytic in a neighborhood of $r^{1/2}$.
Thus
\begin{equation}\label{eq:MwithD}
[x^{2m}]\,f_{k,Q}(x^2) ~=~
D\frac{(mD_Q)^k}{k!}r^{-m} + o(M^kr^{-n})
 ~~\mbox{for some $D,D_Q$.}
\end{equation}

\noindent{\bf The constants in (\ref{eq:MwithD})}.
Since $M(m,Q^0)$ counts compositions with no marked parts, it equals $C(m)$,
the number of compositions of $m$.
Since $C(m)\sim Ar^{-m}$, we have $D=A$.

Note that $M(m,Q^1)/C(m)$ is the average number of parts of size $Q$ in a composition of $m$.
We use the notation of Theorem~\ref{thm:main}(a), the results in Lemma~\ref{lemma:geometric},
and the fact~\cite{BC2} that $X_0(m)$ is asymptotically normal with mean and variance
asymptotically proportional to $m$.
By Lemma~\ref{lemma:geometric}, if we let $P$ (and hence $Q$) go to infinity sufficiently
slowly with $m$, then $\Ex(X_Q(m))/\Ex(X_0(m))\sim Br^Q$.
In this case
$$
mD_Q ~\sim~ \frac{M(m,Q^1)}{M(m,Q^0)} ~=~ \Ex(X_Q(m))
~=~ \frac{\Ex(X_Q(m))}{\Ex(X_0(m))} \Ex(X_0(m))
~\sim~ Br^Q \Ex(X_0(m)).
$$
Recalling the defining relationship
$$
C = B \lim_{m\rightarrow\infty} \frac{\Ex(X_0(m))}{m}
$$
in Theorem~\ref{thm:main}(b),
we have $D_Q\sim Cr^Q$ as $Q\to\infty$.

\medskip

It follows from (\ref{eq:Mdiff}) that, if $P\to\infty$ sufficiently slowly with $m$,
then (\ref{eq:M*M}) holds and so (\ref{eq:MwithD}) provides the asymptotics for
$M(m,Q^k)$ and $M(n,\L)$.
\EndProof

\section{Proofs of Lemmas~\ref{GW12} and \ref{lemma:zeta sum}}

\BeginProofOf{Lemma~\ref{GW12}}
This is Lemma 12 of \cite{GW}.
\EndProof

\BeginProofOf{Lemma~\ref{lemma:zeta sum}}
Fix $\cvec r(x)$ as in Definition~\ref{def:asymp free},
let $P=\max(r_i)$ and $Q\ge\free(P)$.
By Lemma~\ref{lemma:normality}(a),
$\Pr(\zeta_j<k)=o(1)$ for every fixed recurrent $j$.
Thus
$$
\sum_{j<Q}\Pr(\zeta_j<k) ~=~ o(1).
$$

We now consider $j\ge Q$.
Let $\NC(n)$ be the number of compositions of $n$ and
let $\NC_j(n)$ be the number of those having fewer than
$k$ copies of the part $j$.
For some $\delta>0$ to be specified later, let
$\NC^+_j(n)$ be the number of those containing at least $\delta n$
copies of $\cvec r(x)$ and $\NC^-_j(n)$ be the remainder.

Let $\NC^-(n)$ count compositions with fewer than $\delta n$
copies of $\cvec r(x)$.
By Lemma~\ref{lemma:normality}(b) with $\delta=C_1$ sufficiently
small, $\NC^-(n)/\NC(n)$ goes to zero exponentially as $n\to\infty$.
Since $\NC_j^-(n)\le \NC^-(n)$, it follows that
$\NC_j^-(n)< C_2(1+C_3)^{-n}\NC(n)$ where the constants do not
depend on $j$.

In each composition counted by $\NC^+_j(n)$ replace $x$ by
$j$ in $k$ of the $\cvec r(x)$.
This can be done in at least ${\delta n\choose k}$ ways,
giving a composition of $n+k(j-x)$.
Since the resulting composition can have at most $2k-1$ parts of
size $j$, it could have arisen by this replacement process in at
most ${2k-1\choose k}$ ways.
Thus
$$
\NC^+_j(n){\delta n\choose k} ~\le~ \NC(n+k(j-x)){2k-1\choose k}
$$
and so
$$
\frac{\NC^+_j(n)}{\NC(n)}
~\le~ \frac{{2k-1\choose k}}{{\delta n\choose k}}\;
\frac{\NC(n+k(j-x))}{\NC(n)}
< \frac{Br^{-kj}}{n^k}
$$
for some $B=B(k,\delta,x)$ independent of $j$.
Thus
$$
\NC_j(n)/\NC(n) ~<~ Br^{-kj}/n^k + C_2(1+C_3)^{-n},
$$
where all constants are independent of $n$ and $j$.
Summing the right side over $Q\le j\le\log n-\omega(n)$,
we obtain the bound
$$
C_4r^{k\omega(n)} + C_2(1+C_3)^{-n}\log n ~=~ o(1)
$$
for some constants $C_i$.
\EndProof

\section{Proofs of Theorems \ref{thm:AeqC} and \ref{thm:Poisson}}

\BeginProofOf{Theorem~\ref{thm:AeqC}}
There is an increasing function $A_\ell(k,P)$ with supremum $A$
such that the number of compositions of $n\ge k$ with end parts
at most $P$ is at least $A_\ell(k,P)r^{-n}$.
There is a decreasing function $A_u(k)$ with infimum $A$ such that the
number of compositions of $n\ge k$ is at most $A_u(k)r^{-n}$.

Using the construction $\cvec c=\cvec ax\cvec b$ in the statement
of Theorem~\ref{thm:AeqC} together with the idea and notation
in the above proof of Lemma~\ref{lemma:L count},
we construct a composition with one marked part.
If $k\le t\le n-k-Q$, the number of compositions $\cvec aQ\cvec b$
with $\partsum(\cvec a)=t$ and $\partsum(\cvec b)=n-t-Q$
is between $A_\ell(k,P)^2r^{Q-n}$ and $A_u(k)^2r^{Q-n}$.
Sum over all $t$ in the interval $k\le t\le n-k-Q$.
Let $k\to\infty$ sufficiently slowly with $n$.
This shows that $C_P$ in the proof of Lemma~\ref{lemma:L count}
satisfies
$$
A^2 ~=~
\left(\lim_{k\to\infty}A_\ell(k,P)\right)^2
\le~ C_PA ~\le~ \left(\lim_{k\to\infty}A_u(k)\right)^2 =~A^2.
$$
The theorem follows.
\EndProof

\BeginProofOf{Theorem~\ref{thm:Poisson}}
We will show that the three hypotheses (i)-(iii) of Lemma~\ref{GW12}
are satisfied with the choices $\gamma(n)=\log(Cn)$, $\alpha=r$, and
$c=r$.  The first, (i), is obvious.

For (ii), let $\ell,m_1,\dots,m_{\ell}$ be fixed and let
$k_1(n)<k_2(n)<\cdots<k_{\ell}(n)$ be sequences satisfying
$k_i(n)=\log n + O(1)$.  The expectation
$\Ex\prod [\zeta_{k_i}]_{m_i}$, when multiplied by $C(n)$, equals
the number of compositions in which $m_i$ parts of size $k_i$
have been marked and linearly ordered, $1\le i\le\ell$.  Let
$m=\sum_i m_i$ and let
$(L_1,\dots,L_m)$ be one of the ${m \choose m_1,\dots,m_{\ell}}$
possible linear orders of $m_1$ $k_1$'s, etc.  Given a marked
composition counted by Lemma~\ref{lemma:L count}, the linear orders
may be imposed on the marked parts in $\prod m_i!$ ways.  Hence,
$$
C(n) \, \Ex\prod [\zeta_{k_i}]_{m_i}
~\sim~
{m \choose m_1,\dots,m_{\ell}} \, \prod m_i! \,
\frac{A(Cn)^m r^{s-n}}{m!}
$$
where $s=\sum_i m_i k_i$.  Dividing both sides by $C(n)\sim Ar^{-n}$
and noting $Cn=r^{-\gamma(n)}$
completes the confirmation of hypothesis (ii).

Finally, the third hypothesis (iii) is given by
Lemma~\ref{lemma:normality} (c).
\EndProof

\section{Proof of Theorem~\ref{thm:main}}

We recall that all logarithms are to the base $1/r$.

\medskip

\BeginProofOf{Theorem~\ref{thm:main}(a)}
Assertion (a) was proved in Lemma~\ref{lemma:geometric},
except for the formula relating $B$ and $C$.
That relation was proved in the last part of the proof of
Lemma~\ref{lemma:L count}.
\EndProof

\BeginProofOf{Theorem~\ref{thm:main}(d)}
The final claim in (d) is easily proved by bounding the right side
of~(\ref{eq:qn(k)}).
One can also show that the sum of the right side of~(\ref{eq:qn(k)})
over $|k-\log n|>\omega_d(n)$ tends to zero.
Thus it suffices to prove~(\ref{eq:qn(k)}).

By Theorem~\ref{thm:Poisson} there is some $\omega(n)\to\infty$
such that
\begin{eqnarray*}
q_n(k)
&\sim& p(n)
\Biggl(\prod_{j=f(n)}^k\left(1-\exp\left(-r^{j-\log(Cn)}\right)\right)\Biggr)
\Biggl(\prod_{j=k+1}^n\exp\left(-r^{j-\log(Cn)}\right)\Biggr)\\
&=& p(n)
\Biggl(\prod_{j=f(n)}^k\left(1-\exp\left(-Cnr^j\right)\right)\Biggr)
\exp\Biggl(-Cn\sum_{j=k+1}^n r^j\Biggr),
\end{eqnarray*}
where $f(n)=\lfloor\log(Cn)-\omega(n)\rfloor$ and $p(n)$ is the probability
that a random composition of~$n$ contains all recurrent parts less than $f(n)$.
With a little calculation, we see that Theorem~\ref{thm:main}(d)
is equivalent to $p(n)\sim1$.
With $\zeta_j$ as in Lemma~\ref{lemma:zeta sum}, we have
$$
p(n) ~\ge~ 1-\sum_{j\le f(n)}\Pr(\zeta_j\!=\!0) ~=~ 1-o(1).
$$
We now prove the final claim about strong concentration.
Since the probability of being gap-free is bounded away from zero
by previous claim in (d), it suffices to prove the validity of the
statement for all compositions. The result now follows from (b).
This proves Theorem~\ref{thm:main}(d).
\EndProof

\BeginProofOf{Theorem~\ref{thm:main}(b,c)}
We first turn to the formula in (c).
Since $\zeta_j$ in Theorem~\ref{thm:Poisson} is irrelevant for small $j$,
we let $\zeta_j$ be as in Lemma~\ref{lemma:zeta sum} and
$\omega(n)$ be any function that goes to infinity.
Note that
$$
\Ex(D_n)+\nu ~=~ \sum_{j=1}^n \Pr(\zeta_j\ne 0)
$$
and so, by Lemma~\ref{lemma:zeta sum},
\begin{eqnarray}
\Ex(D_n)+\nu
&=& o(1) ~+\!\! \sum_{j=1}^{\kpn-1}\!\!1
\label{eq:Dn split 1}\\
&& +\!\! \sum_{j=\kmn}^{\kpn}\!\!\Pr(\zeta_j\ne 0)
\label{eq:Dn split 2}\\
&& +\!\! \sum_{j=\kpn+1}^n\!\!\Pr(\zeta_j\ne 0).
\label{eq:Dn split 3}
\end{eqnarray}
By Lemma~\ref{lemma:normality}(c)
$$
\Pr(\zeta_j\ne0) ~=~ O(nr^j)
~~\mbox{provided $j\to\infty$ with $n$.}
$$
Thus the sum in (\ref{eq:Dn split 3}) is
$O(r^{\omega(n)})=o(1)$.
By the Poisson distribution, the terms in the sum (\ref{eq:Dn split 2})
are $1-\exp\left(-Cnr^j\right)+o(1)$ and so, if $\omega(n)\to\infty$
sufficiently slowly, that sum is
$$
o(1)+\sum\Bigl(1-\exp\left(-Cnr^j\right)\Bigr).
$$
Since the sum of $\exp\left(-Cnr^j\right)$ over $j<\kmn$ is $o(1)$,
we may replace 1 in the sum (\ref{eq:Dn split 1}) with
$1-\exp\left(-Cnr^j\right)$.
Finally,
$$
\sum_{j>\kpn}\!\!\Bigl(1-\exp\left(-Cnr^j\right)\Bigr) ~=~ o(1)
$$
and so
$$
\Ex(D_n)+\nu
~=\! \sum_{j\ge 0}\Bigl(1-\exp\left(-Cnr^j\right)\Bigr)-1+o(1).
$$
Let  $f(x) = \sum_{j\ge 0}\Bigl(1-\exp\left(-xr^j\right)\Bigr)$.
Then $\Ex(D_n)=f(Cn)-1+o(1)$.
We use the standard Mellin transform.
(See \cite[p.765]{FS}, and also their Example B.5 which
treats $r=1/2$).
It follows that
$$
f(x) ~=~ \log x+\gamma\log e+\frac{1}{2}+P_0(x)+o(1),
$$
where $P_0(x)$ is given
by~(\ref{eq:oscillation}).
This proves Theorem~\ref{thm:main}(c).

For the maximum part size $M_n$, we proceed in a similar manner:
\begin{eqnarray*}
\Ex(M_n)
&=&\sum_{j=1}^n\Pr(M_n\ge j)\\
&=&\!\sum_{j=1}^{\kmn-1}\!1
~+\!\!\! \sum_{j=\kmn}^{n}\!\Bigl(1-\Pr(\wedge_{i\ge j} \{\zeta_j=0\})\Bigr)+o(1)\\
&=&\sum_{j=1}^{\kmn-1}\hskip-12pt 1 ~~+\!\!
\sum_{j=\kmn}^{n}\!\!\left(1-\exp\left(-\frac{Cn}{1-r}r^j\right)\right)+o(1)\\
&=&\sum_{j\ge 0}^{n}\left(1-\exp\left(-\frac{Cn}{1-r}r^j\right)\right)-1+o(1)
=f\left(\frac{Cn}{1-r}\right)-1+o(1)\\
&=& \log\left(\frac{Cn}{1-r}\right) + \gamma\log e
- \frac{1}{2}+P_0\left(\frac{Cn}{1-r}\right)+o(1).
\end{eqnarray*}

It remains to prove the claims about $|D_n-\log n|$ and $|M_n-\log n|$.
Since $D_n\le M_n$, it suffices to establish the lower bound for $D_n$ and
the upper for $M_n$.
The upper bound on $M_n$ was proved in \cite[Section~9]{BC2}.
In the proof of Theorem~1(d) we showed that $p(n)\sim1$,
which establishes the lower bound on $D_n$.
\EndProof

\BeginProofOf{Theorem~\ref{thm:main}(e--g)}
Let $\Gamma:=(\Gamma_1,\Gamma_2,\ldots,\Gamma_m)$ be a
sequence of i.i.d.\ geometric random variables with parameter
$p=1-r$.
Hitczenko and Knopfmacher~\cite{HK} showed that the probability
the sequence $\Gamma$ is gap-free is given by the $p_m$
in our~(\ref{eq:pm}) and they established the oscillation of $p_m$
when $p\ne1/2$.

Let $\omega(m)$ go to infinity arbitrarily slowly with $m$.
Let $M_m$ be the largest $\Gamma_i$.

By Theorem~\ref{thm:main}(b) $|M_m-\log m|<\omega(m)$, and,
as was shown in the proof of Theorem~1(d), all
recurrent parts less than $\log m-\omega(m)$ are
asymptotically almost surely present in
$\Gamma$.  Let
$$
\matrix{ \zeta_j:=|\{i:\Gamma_i=j\}|_{\vphantom{\bigm|}},\quad&
\lambda_j:=m(1-r)r^{j-1},\hfill \cr
 k^-:=\km,\quad& k^+:=\kp.\cr}
$$
When $k^-\le k\le k^+$,
$$
\Pr(\zeta_j=k) ~\sim~ e^{-\lambda_j}\lambda_j^k/k!
$$
by the standard Poisson approximation for i.i.d.\ rare random
variables.
It should be well-known that $\{\zeta_j: k^-\le j\le k^+\}$
are asymptotically independent, but we include a proof since
we lack a reference.
For all fixed positive integers $m_1,\ldots, m_j$, we have
\begin{eqnarray*}
\Pr\left(\wedge_{k^-\le j\le k^+}\{\zeta_j=m_j\}\right)
&=&\frac{m! (\lambda_{k^-}/m)^{m_{k^-}}\cdots(\lambda_{k^+}/m)^{m_{k^+}}}
{(m_{k^-})!\cdots (m_{k^+})!(m-m_{k^-}+\cdots+m_{k^+})!}\\
& &\times
\left(1-\frac{\lambda_{k^-}+\ldots+\lambda_{k^+}}{m}\right)^{m-(m_{k^-}
+\ldots+m_{k^+})}\\
&\sim & \frac{\lambda_{k^-}^{m_{k^-}}\cdots
\lambda_{k^+}^{m_{k^+}}}{(m_{k^-})!\cdots
(m_{k^+})!}\exp(-(\lambda_{k^-}+\ldots +\lambda_{k^+})).
\end{eqnarray*}

 Thus, with $k$ the largest part,
\begin{eqnarray}
\nonumber p_m &\sim& \sum_{k=k^-}^{k^+}
\Biggl(\prod_{j=k+1}^{k^+} e^{-\lambda_j}\Biggr)
\Biggl(\prod_{j=k^-}^k (1-e^{-\lambda_j})\Biggr)\\
&\sim& \label{eq:pm new}
\sum_{k=k^-}^{k^+} \exp\left(-mr^k\right)
\prod_{j=k^-}^k\left(1-\exp\left(-m(1-r)r^{j-1}\right)\right).
\end{eqnarray}

Equation~(\ref{eq:pm new}) is the same as the sum of (\ref{eq:qn(k)})
if $m=Cn/(1-r)$.
However, (\ref{eq:pm new}) was derived under the assumption that $m$
is an integer.
We now treat (\ref{eq:pm new}) as a function of real variable $m$,
say $f(m)$, and show that $f'(m)=o(1)$ as $m\to\infty$.
It then follows that $f(x)\sim f(\lfloor x\rfloor)$ as $x\to\infty$
and we will be done.
Call the terms in the sum (\ref{eq:pm new}) $T_k(m)$.
We have
\begin{eqnarray*}
|T_k'(m)|
&<& \left|\frac{T_k'(m)}{T_k(m)}\right|
~=~ \left|\left(\ln(T_k(m)\right)'\right|\le r^k
+ \sum_{j=k^-}^k\frac{(1-r)r^{j-1}}{\exp\left(m(1-r)r^{j-1}\right)-1} \\
&<& r^k +  \sum_{j=k^-}^k\frac{(1-r)r^{j-1}}{m(1-r)r^{j-1}}
\le r^{k^-} + \frac{k-k^-+1}{m}
~<~ \frac{\omega_1(m)}{m}
\end{eqnarray*}
for some $\omega_1(m)\to\infty$ much slower than $m$.
Since there are only $2\omega(m)$ values for $k$, $f'(m)=o(1)$.

The oscillation is associated with the imaginary poles, which are at
$2k\pi i/\ln(1/q)$ in the notation of~\cite{HK}.
When the result is translated back from $m$ to $n$, we obtain
the same period as $P$ in~(\ref{eq:oscillation}).

\medskip

We now prove (f).
It follows from Theorem~\ref{thm:Poisson} that
$$
g_n(k)~\sim\hskip-.05in
\sum_{j>\log(Cn)-\omega(n)}\hskip-.1in
\Pr(\zeta_{j}=k,\zeta_{j+1}=\zeta_{j+2}=\cdots=0).
$$
Setting $j=\ell+\lfloor\log(Cn)\rfloor$ and
$\delta(n)=Cn-\lfloor\log(Cn)\rfloor$,
$$
g_n(k)\sim
\sum_{\ell=-\infty}^{\infty}\frac{r^{k(\ell-\delta(n))}}{k!}
\prod_{i\ge\ell}\exp\left(-r^{i-\delta(n)}\right)
\sim
\sum_{\ell=-\infty}^{\infty}\frac{r^{k(\ell-\delta(n))}}{k!}
\exp\left(\frac{-r^{\ell-\delta(n)}}{1-r}\right).
$$
It follows from Poisson's summation formula \cite{Olver} that
$$
g_n(k) ~\sim~
\sum_{\ell=-\infty}^{\infty}\int_{-\infty}^{\infty}
\frac{1}{k!}\exp(-2\pi i\ell t)r^{k(t-\delta(n))}
\exp\left(\frac{-r^{t-\delta(n)}}{1-r}\right)dt.
$$
Setting $z=\frac{r^{t-\delta(n)}}{1-r}$,
\begin{eqnarray*}
g_n(k)&\sim&
\frac{(1-r)^k\log e}{k!}\sum_{\ell=-\infty}^{\infty}
  \exp(-2\pi i\ell(\delta(n)-\log(1-r)))
  \int_{0}^{\infty}e^{-z}z^{k-1+2\pi i\ell\log e}dz\\
&\sim & \frac{(1-r)^k\log e}{k!}\sum_{\ell=-\infty}^{\infty}
  \Gamma\left(k+2\pi i\ell\log e\right)
  \exp\left(-2\pi i\ell\log\frac{Cn}{1-r}\right)\\
&\sim &
\frac{(1-r)^k}{k!}P_k\left(\frac{Cn}{1-r}\right)+\frac{(1-r)^k\log
e}{k}
\end{eqnarray*}
This completes the proof of (f).

\medskip

We now prove (g).
By Lemma~\ref{lemma:zeta sum} and Theorem~\ref{thm:main}(b)
we may limit our attention to parts $j$ for which
$|j-\log n|\le \omega(n)$.
By Theorem~\ref{thm:Poisson}, the probability that part $j$
appears with multiplicity $k$ is asymptotically
$e^{-\mu_j}\mu_j^k/k!$ where $\mu_j=Cnr^j$.
Using the Poisson summation formula as in the proof of (f),
the expected number of parts of multiplicity $k$ is asymptotic to
\begin{eqnarray*}
& &\frac{1}{k!}\sum_j \exp\left(-Cnr^j\right)(Cnr^j)^k\\
&\sim
&\frac{1}{k!}\sum_{\ell=-\infty}^{\infty}\int_{-\infty}^{\infty}\exp\left(-2i\pi
\ell t-r^{t-\delta(n)}\right)r^{k(t-\delta(n))}dt\\
&\sim &\frac{\log e}{k!}\sum_{\ell=-\infty}^{\infty}\exp\left(-2i\pi
\ell \log(Cn)\right)\Gamma\left(k+2i\pi \ell\log e\right)\\
&\sim &\frac{P_k(Cn)}{k!}+\frac{\log e}{k}.
\end{eqnarray*}
The claim about $m_n(k)$ follows from the fact that
$m_n(k)=\Ex(D_n(k)/D_n)$ and the tight concentration of $D_n$ in (c)---an
argument used by Louchard~\cite{Lo} for unrestricted compositions.
\EndProof

\vskip 20pt

{}

\end{document}